\documentclass[10pt,leqno]{article}

\usepackage{amsmath,amssymb,amsthm,mathrsfs,dsfont}

\usepackage[margin=3cm]{geometry} 

\usepackage{titlesec,hyperref}

\usepackage{color}

\usepackage{fancyhdr}
\titleformat{\subsection}{\it}{\thesubsection.\enspace}{1pt}{}

\newtheorem{theo}{Theorem}[section]
\newtheorem{lemm}[theo]{Lemma}

\newtheorem{rema}[theo]{Remark}
\numberwithin{equation}{section}

\allowdisplaybreaks 

\newcommand\ep{{\varepsilon}} 
\newcommand\lm{{\lesssim}}

\begin{document}
\title{The $L^2$ decay for the 2D co-rotation FENE dumbbell model of polymeric flows
\hspace{-4mm}
}

\author{Wei $\mbox{Luo}^1$\footnote{E-mail:  luowei23@mail2.sysu.edu.cn} \quad and\quad
 Zhaoyang $\mbox{Yin}^{1,2}$\footnote{E-mail: mcsyzy@mail.sysu.edu.cn}\\
 $^1\mbox{Department}$ of Mathematics,
Sun Yat-sen University, Guangzhou 510275, China\\
$^2\mbox{Faculty}$ of Information Technology,\\ Macau University of Science and Technology, Macau, China}

\date{}
\maketitle
\hrule

\begin{abstract}
In this paper we mainly study the long time behaviour of solutions to the finite extensible nonlinear elastic (FENE) dumbbell model with dimension two in the co-rotation case. Firstly, we obtain the $L^2$ decay rate of the velocity of the 2D co-rotation FENE model is $(1+t)^{-\frac{1}{2}}$ with small data. Then, by virtue of the Littlewood-Paley theory, we can remove the small condition. Our obtained sharp result improves considerably the recent results in \cite{Luo-Yin,Schonbek}. \\

\vspace*{5pt}

\vspace*{5pt}
\end{abstract}

\vspace*{10pt}

\tableofcontents

\section{Introduction}
   In this paper we consider the finite extensible nonlinear elastic (FENE) dumbbell model \cite{Bird1977,Doi1988}:
   \begin{align}
\left\{
\begin{array}{ll}
u_t+(u\cdot\nabla)u-\nu\Delta{u}+\nabla{P}=div~\tau, ~~~~~~~div~u=0,\\[1ex]
\psi_t+(u\cdot\nabla)\psi=div_{R}[- \sigma(u)\cdot{R}\psi+\beta\nabla_{R}\psi+\nabla_{R}\mathcal{U}\psi],  \\[1ex]
\tau_{ij}=\int_{B}(R_{i}\nabla_{j}\mathcal{U})\psi dR, \\[1ex]
u|_{t=0}=u_0,~~\psi|_{t=0}=\psi_0, \\[1ex]
(\beta\nabla_{R}\psi+\nabla_{R}\mathcal{U}\psi)\cdot{n}=0 ~~~~ \text{on} ~~~~ \partial B(0,R_{0}) .\\[1ex]
\end{array}
\right.
\end{align}
In (1.1)~~$\psi(t,x,R)$ denotes the distribution function for the internal configuration and $u(t,x)$ stands for the velocity of the polymeric liquid, where $x\in\mathbb{R}^{d}$ and $d\geq2$ means the dimension. Here the polymer elongation $R$ is bounded in ball $ B=B(0,R_{0})$ of $\mathbb{R}^{d}$ which means that the extensibility of the polymers is finite. $\beta=\frac{2k_BT_a}{\lambda}$, where $k_B$ is the Boltzmann constant, $T_a$ is the absolute temperature and $\lambda$ is the friction coefficient. $\nu>0$ is the viscosity of the fluid, $\tau$ is an additional stress tensor and $P$ is the pressure. The Reynolds number $Re=\frac{\gamma}{\nu}$ with $\gamma\in(0,1)$ and the density $\rho=\int_B\psi dR$. Moreover the potential $\mathcal{U}(R)=-k\log(1-(\frac{|R|}{|R_{0}|})^{2})$ for some constant $k>0$. $\sigma(u)$ is the drag term. In general, $\sigma(u)=\nabla u$. For the co-rotation case, $\sigma(u)=\frac{\nabla u-(\nabla u)^{T}}{2}$.

   This model describes the system coupling fluids and polymers. The system is of great interest in many branches of physics, chemistry, and biology, see \cite{Bird1977,Doi1988}. In this model, a polymer is idealized as an "elastic dumbbell" consisting of two "beads" joined by a spring that can be modeled by a vector $R$. At the level of liquid, the system couples the Navier-Stokes equation for the fluid velocity with a Fokker-Planck equation describing the evolution of the polymer density. This is a micro-macro model (For more details, one can refer to  $\cite{Bird1977}$, $\cite{Doi1988}$, $\cite{Masmoudi2008}$ and $\cite{Masmoudi2013}$).

%
%
%


In the paper we will take $\beta=1$, $\nu=1$ and $R_{0}=1$.
Notice that $(u,\psi)$ with $u=0$ and $$\psi_{\infty}(R)=\frac{e^{-\mathcal{U}(R)}}{\int_{B}e^{-\mathcal{U}(R)}dR}=\frac{(1-|R|^2)^k}{\int_{B}(1-|R|^2)^kdR},$$
is a trivial solution of (1.1). By a simple calculation, we can rewrite (1.1) for the following system:
\begin{align}
\left\{
\begin{array}{ll}
u_t+(u\cdot\nabla)u-\Delta u+\nabla{P}=div~\tau,  ~~~~~~~div~u=0,\\[1ex]
\psi_t+(u\cdot\nabla)\psi=div_{R}[-\sigma(u)\cdot{R}\psi+\psi_{\infty}\nabla_{R}\frac{\psi}{\psi_{\infty}}],  \\[1ex]
\tau_{ij}=\int_{B}(R_{i}\nabla_{R_j}\mathcal{U})\psi dR, \\[1ex]
u|_{t=0}=u_0, \psi|_{t=0}=\psi_0, \\[1ex]
\psi_{\infty}\nabla_{R}\frac{\psi}{\psi_{\infty}}\cdot{n}=0 ~~~~ \text{on} ~~~~ \partial B(0,1) .\\[1ex]
\end{array}
\right.
\end{align}
{\bf Remark.} As in the reference \cite{Masmoudi2013}, one can deduce that $\psi=0$ on the boundary.

There are a lot of mathematical results about the FENE dumbbell model. M. Renardy \cite{Renardy} established the local well-posedness in Sobolev spaces with potential $\mathcal{U}(R)=(1-|R|^2)^{1-\sigma}$ for $\sigma>1$. Later, B. Jourdain, T. Leli\`{e}vre, and
C. Le Bris \cite{Jourdain} proved local existence of a stochastic differential equation with potential $\mathcal{U}(R)=-k\log(1-|R|^{2})$ in the case $k>3$ for a Couette flow. H. Zhang and P. Zhang \cite{Zhang-H} proved local well-posedness of (1.4) with $d=3$ in weighted Sobolev spaces. For the co-rotation case, F. Lin, P. Zhang, and Z. Zhang \cite{F.Lin} obtain a global existence results with $d=2$ and $k > 6$. If the initial data is perturbation around equilibrium, N. Masmoudi \cite{Masmoudi2008} proved global well-posedness of (1.4) for $k>0$. In the co-rotation case with $d=2$, he \cite{Masmoudi2008} obtained a global result for $k>0$ without any small conditions. In the co-rotation case, A. V. Busuioc, I. S. Ciuperca, D. Iftimie and L. I. Palade \cite{Busuioc} obtain a global existence result with only the small condition on $\psi_0$. The global existence of weak solutions in $L^2$ was proved recently by N. Masmoudi \cite{Masmoudi2013} under some entropy conditions.

Recently, M. Schonbek \cite{Schonbek} studied the $L^2$ decay of the velocity for the co-rotation
FENE dumbbell model, and obtained the
decay rate $(1+t)^{-\frac{d}{4}+\frac{1}{2}}$, $d\geq 2$ with $u_0\in L^1$.
Moreover, she conjectured that the sharp decay rate should be $(1+t)^{-\frac{d}{4}}$,~$d\geq 2$.
However, she failed to get it because she could not use the bootstrap argument as in \cite{Schonbek1985} due to the
additional stress tensor.  More recently, W. Luo and Z. Yin \cite{Luo-Yin} improved Schonbek's result
and showed that the decay rate is $(1+t)^{-\frac{d}{4}}$ with $d\geq 3$ and $\ln^{-l}(1+t)$ with $d=2$ for any $l\in\mathbb{N^+}$.
This result shows that M. Schonbek's conjecture is true when $d \geq3$. However, there is no any result to show that M. Schonbek's conjecture is true when $d=2$.

In this paper, we are going on to prove that M. Schonbek's conjecture holds true when $d=2$.
Firstly, we show that the $L^2$ decay rate is $(1+t)^{-\frac{d}{4}}$, $d\geq 2$ for the velocity with the small initial data.
The main idea is that we can obtain a estimate for the $L^3_T(L^2)-$ norm of the velocity with the initial data being small in $L^2$.
Next, we can use the bootstrap argument to obtain the optimal decay rate.
Since we are interested in the
large time behaviour of the weak solutions,
we can consider the evolution system after a large time $T_0$.
 Note that the $L^2$ norm of the solutions
 will decay to zero. Thus we can remove the small $L^2$ norm condition of the initial data.
 And then the main difficult is to get the $L^1$-estimate
for the velocity. Since the Leray project operator is not bounded in $L^1$, it follows that one can not obtain the $L^1$-estimate directly by
the heat kernel estimate. Instead, we will use the Littlewood-Paley theory to estimate the $\dot{B}^0_{1,1}$ norm in place of the $L^1$ norm. Finally, we can prove that the sharp $L^2$ decay rate is $(1+t)^{-\frac{d}{4}}$, $d\geq 2$ for the velocity with the large initial data.

  The paper is organized as follows. In Section 2 we introduce some notations and  give some preliminaries which will be used in the sequel. In Section 3 we study the $L^2$ decay of solutions to the 2D co-rotation FENE model
  by using the Fourier splitting method and the Littlewood-Paley theory.

\section{Preliminaries}
  In this section we will introduce some notations and useful lemmas which will be used in the sequel.

 If the function spaces are over $\mathbb{R}^d$ and $B$ with respect to the variable $x$ and $R$, for simplicity, we drop $\mathbb{R}^d$ and $B$ in the notation of function spaces if there is no ambiguity.

For $p\geq1$, we denote by $\mathcal{L}^{p}$ the space
$$\mathcal{L}^{p}=\big\{\psi \big|\|\psi\|^{p}_{\mathcal{L}^{p}}=\int \psi_{\infty}|\frac{\psi}{\psi_{\infty}}|^{p}dR<\infty\big\}.$$

  We will use the notation $L^{p}_{x}(\mathcal{L}^{q})$ to denote $L^{p}[\mathbb{R}^{d};\mathcal{L}^{q}]:$
$$L^{p}_{x}(\mathcal{L}^{q})=\big\{\psi \big|\|\psi\|_{L^{p}_{x}(\mathcal{L}^{q})}=(\int_{\mathbb{R}^{d}}(\int_{B} \psi_{\infty}|\frac{\psi}{\psi_{\infty}}|^{q}dR)^{\frac{p}{q}}dx)^{\frac{1}{p}}<\infty\big\}.$$
When $p=q$, we also use the short notation $\mathcal{L}^p$ for $L^p_x(\mathcal{L}^{p})$ if there is no ambiguity.

The symbol $\widehat{f}=\mathcal{F}(f)$ denotes the Fourier transform of $f$.

Moreover, we denote by $\dot{\mathcal{H}}^1$ the space
$$\dot{\mathcal{H}}^1=\big\{g\big| \|g\|_{\dot{\mathcal{H}}^1}=(\int_B|\nabla_R g|^2\psi_\infty dR)^{\frac{1}{2}}\big\}.$$
Sometimes we write $f\lm g$ instead of $f\leq Cg$, where $C$ is a constant. We agree that $\nabla$ stands for $\nabla_x$ and $div$ stands for $div_x$.

The following lemma allows us to estimate the extra stress tensor $\tau$.
\begin{lemm}\cite{Masmoudi2008}\label{Lemma2}
 If $\int_B \psi dR=0$ and
   $\displaystyle\int_B\bigg|\nabla _R  \bigg(\displaystyle\frac{\psi}{\psi _\infty }\bigg)\bigg|^2{\psi _\infty} dR<\infty$ with $p\geq 2$, then there exists a constant $C$ such that
   \[\int_{B}\frac{|\psi|^{2}}{\psi_{\infty}}dR\leq C \displaystyle\int_B\bigg|\nabla _R  \bigg(\displaystyle\frac{\psi}{\psi _\infty }\bigg)\bigg|^2{\psi _\infty} dR.\]
\end{lemm}
\begin{lemm}\label{Lemma3}
\cite{Masmoudi2008} For all $\varepsilon>0$, there exists a constant $C_{\varepsilon}$ such that
$$|\tau|^2\leq\varepsilon\int_{B}\psi_{\infty}|\nabla_{R}\frac{\psi}{\psi_{\infty}}|^{2}dR
+C_{\varepsilon}\int_{B}\frac{|\psi|^{2}}{\psi_{\infty}}dR,$$
or
$$|\tau|^2\leq C(\int_{B}\frac{|\psi|^{2}}{\psi_{\infty}}dR)^{\frac{1}{2}}(\int_{B}\psi_{\infty}|\nabla_{R}\frac{\psi}{\psi_{\infty}}|^{2}dR)^{\frac{1}{2}}.$$
\end{lemm}

The following lemma is the well-known $(L^p,L^q)-$estimates which can be easily deduced from the properties of the heat kernel.
\begin{lemm}\cite{Bahouri2011}\label{Lemma4}
Let $1\leq p\leq q\leq \infty$. For all $f\in L^p$, there exists a constant $C$ such that
$$\|e^{-t\Delta}f\|_{L^q}\leq Ct^{-\frac{d}{2}(\frac{1}{p}-\frac{1}{q})}\|f\|_{L^p}, \quad \|e^{-t\Delta}\nabla f\|_{L^q}\leq Ct^{-\frac{d}{2}(\frac{1}{p}-\frac{1}{q})-\frac{1}{2}}\|f\|_{L^p}.$$
\end{lemm}

\begin{theo}\cite{Luo-Yin}\label{th1}
Let $(u,\psi)$ be a weak solution of (1.2) with the initial data $u_0\in L^2\cap L^1$ and $\psi_0$ satisfies $\psi_0-\psi_\infty\in L^2_x(\mathcal{L}^2)$ and $\int_B\psi_0=1$ $a.e.$ in $x$. Then there exists a constant $C$ such that
\begin{align}
\int_{\mathbb{R}^d\times B}\frac{|\psi-\psi_\infty|^2}{\psi_\infty}dxdR\leq C\exp{(-Ct)},
\end{align}
\begin{align}
\|u\|_{L^2}\leq C(1+t)^{-\frac{d}{4}},\quad \text{if} \quad d\geq3, \quad \|u\|_{L^2}\leq C_l\ln^{-l}(e+t),\quad \text{if}\quad d=2,
\end{align}
where $l>0$ is arbitrarily integer and $C_l$ is a constant dependent on $l$.
\end{theo}

\section{Main results}
This section is devoted to investigating the long time behaviour for the velocity of the co-rotation FENE dumbbell model with dimension $d=2$. More precisely, we prove the $L^2$ decay for the weak solutions of the 2D co-rotation FENE dumbbell model and obtain the $L^2$ decay rate. The existence of the solutions in $L^2$ was established in \cite{Lions-Masmoudi2007,Schonbek}. Then our main result can be stated as follows.
\subsection{The $L^2$ decay with small data}
\begin{theo}\label{th2}
Let $(u,\psi)$ be a weak solution of (1.2) with the initial data $u_0\in L^2\cap L^1$ and $\psi_0$ satisfies $\psi_0-\psi_\infty\in L^2_x(\mathcal{L}^2)$ and $\int_B\psi_0=1$ $a.e.$ in $x$. A constant $\ep$ exists such that if
$$\|u_0\|_{L^2}+\|\psi_0-\psi_\infty\|_{\mathcal{L}^2}\leq \ep, $$
then we have
\begin{align}
\int_{\mathbb{R}^2\times B}\frac{|\psi-\psi_\infty|^2}{\psi_\infty}dxdR\leq C\exp{(-Ct)},
\end{align}
\begin{align}
\|u\|_{L^2(\mathbb{R}^2)}\leq C(1+t)^{-\frac{1}{2}},
\end{align}
where $C$ is a constant dependent on the initial data.
\begin{proof}
By the standard density argument, we only need to prove that the conclusion holds for the smooth solution.
Since $\psi_\infty=\displaystyle\frac{(1-|R|^2)^k}{\int_{B}(1-|R|^2)^kdR}=\displaystyle\frac{(1-|R|^2)^k}{C_0}$, it follows that
\begin{multline}\label{3.3}
div_R([(\nabla u-(\nabla u)^T]R\psi_\infty)=\sum_{i,j}\partial_{R_i}[(\partial_iu^j-\partial_ju^i)R_j\psi_\infty]
\\
=\sum_{i,j}(\partial_iu^j-\partial_ju^i)\delta_{ij}\psi_\infty+\sum_{i,j}\frac{2k(\partial_iu^j-\partial_ju^i)R_jR_i(1-|R|^2)^{k-1}}{C_0}
=0.
\end{multline}
By virtue of the second equation of (1.2), we have
\begin{equation}\label{3.4}
(\psi-\psi_\infty)_t+(u\cdot\nabla)(\psi-\psi_\infty)=div_{R}[-\sigma(u)\cdot{R}(\psi-\psi_\infty)+\psi_{\infty}\nabla_{R}\frac{\psi-\psi_\infty}{\psi_{\infty}}].
\end{equation}
Multiplying $\frac{\psi-\psi_\infty}{\psi_\infty}$ by both sides of the above equation and integrating over $B$ with $R$, we obtain
\begin{multline}\label{3.5}
\frac{1}{2}\frac{d}{dt}\int_B\frac{|\psi-\psi_\infty|^2}{\psi_\infty}+\frac{1}{2}u\cdot\nabla_x\int_B\frac{|\psi-\psi_\infty|^2}{\psi_\infty}+\int_B\psi_\infty|\nabla_R(\frac{\psi-\psi_\infty}{\psi_\infty})|^2
\\
=\int_B\sigma(u)R(\psi-\psi_\infty)\nabla_R(\frac{\psi-\psi_\infty}{\psi_\infty}).
\end{multline}
Using integration by parts and (\ref{3.3}), we see that
\begin{multline}\label{3.6}
\int_B\sigma(u)R(\psi-\psi_\infty)\nabla_R(\frac{\psi-\psi_\infty}{\psi_\infty})=\int_B\sigma(u)R\psi_\infty[\frac{1}{2}\nabla_R(\frac{\psi-\psi_\infty}{\psi_\infty})^2]\\
=-\frac{1}{2}\int_Bdiv_R([(\nabla u-(\nabla u)^T]R\psi_\infty)(\frac{\psi-\psi_\infty}{\psi_\infty})^2=0.
\end{multline}
Plugging (\ref{3.6}) into (\ref{3.5}) and using the fact that $div~ u=0$, we deduce that
\begin{align}\label{3.7}
\frac{1}{2}\frac{d}{dt}\int_{\mathbb{R}^2\times B}\frac{|\psi-\psi_\infty|^2}{\psi_\infty}+\int_{\mathbb{R}^2\times B}\psi_\infty|\nabla_R(\frac{\psi-\psi_\infty}{\psi_\infty})|^2=0.
\end{align}
By virtue of the equation (1.2), we have $\int_B\psi dR=\int_B\psi_0 dR=1$, which leads to $\int_B(\psi-\psi_\infty) dR=0$. Taking advantage of Lemma \ref{Lemma2}, we infer that
\begin{align}
\frac{1}{2}\frac{d}{dt}\int_{\mathbb{R}^2\times B}\frac{|\psi-\psi_\infty|^2}{\psi_\infty}+C\int_{\mathbb{R}^2\times B}\frac{|\psi-\psi_\infty|^2}{\psi_\infty}\leq 0,
\end{align}
which leads to
\begin{align}
\frac{d}{dt}\bigg[\exp{(Ct)}\int_{\mathbb{R}^2\times B}\frac{|\psi-\psi_\infty|^2}{\psi_\infty}\bigg]\leq 0\Rightarrow \int_{\mathbb{R}^2\times B}\frac{|\psi-\psi_\infty|^2}{\psi_\infty}\leq \exp{(-Ct)}\int_{\mathbb{R}^2\times B}\frac{|\psi_0-\psi_\infty|^2}{\psi_\infty}.
\end{align}
Since $\partial_x \psi_\infty=0$, it follows that $div \tau= div \int_{B}(R\otimes\nabla_{R}\mathcal{U})\psi dR=div \int_{B}(R\otimes\nabla_{R}\mathcal{U})(\psi-\psi_\infty) dR$. Then, we may assume that $\tau =\int_{B}(R\otimes\nabla_{R}\mathcal{U})(\psi-\psi_\infty) dR$ .
By the standard energy estimate for the Navier-Stokes equations, we get
\begin{align}
\frac{1}{2}\frac{d}{dt}\|u\|^2_{L^2}+\|\nabla u\|^2_{L^2}=-\int_{\mathbb{R}^2}\tau:\nabla u\leq \frac{1}{2}\|\nabla u\|_{L^2}+\frac{1}{2}\|\tau\|^2_{L^2}.
\end{align}
Using Lemmas \ref{Lemma2}-\ref{Lemma3}, we verify that
\begin{align}
\frac{d}{dt}\|u\|^2_{L^2}+\|\nabla u\|^2_{L^2}\leq\|\tau\|^2_{L^2}\leq K\int_{\mathbb{R}^2\times B}\psi_\infty|\nabla_R(\frac{\psi-\psi_\infty}{\psi_\infty})|^2.
\end{align}
Let $\lambda\geq 2K$ be a sufficiently large constant. From the above inequality and (\ref{3.7}), we deduce that
\begin{align}\label{3.12}
\frac{d}{dt}(\lambda\|\psi-\psi_\infty\|^2_{\mathcal{L}^2}+\|u\|^2_{L^2})+\lambda\int_{\mathbb{R}^2\times B}\psi_\infty|\nabla_R(\frac{\psi-\psi_\infty}{\psi_\infty})|^2\leq 0.
\end{align}
Taking $\lambda=2K$, we have
\begin{align}
\|u\|^2_{L^2}+\int^t_0\|\nabla u\|^2_{L^2}ds\leq \|u_0\|^2_{L^2}+2K\|\psi_0-\psi_\infty\|^2_{\mathcal{L}^2}<\infty.
\end{align}
From (\ref{3.12}), we have
\begin{multline}
\frac{d}{dt}((1+t)^2\lambda\|\psi-\psi_\infty\|^2_{\mathcal{L}^2}+(1+t)^2\|\widehat{u}\|^2_{L^2})+\lambda (1+t)^2\int_{\mathbb{R}^2\times B}\psi_\infty|\nabla_R(\frac{\psi-\psi_\infty}{\psi_\infty})|^2+(1+t)^2\int_{\mathbb{R}^2}|\xi|^2|\widehat{u}|^2d\xi\\
\leq 2(1+t)\lambda\|\psi-\psi_\infty\|^2_{\mathcal{L}^2}+2(1+t)\|\widehat{u}\|^2_{L^2}.
\end{multline}
 Setting $S(t)=\{\xi:|\xi|^2\leq \frac{2}{1+t}\}$, then we obtain
\begin{multline}\label{3.15}
\frac{d}{dt}((1+t)^2\lambda\|\psi-\psi_\infty\|^2_{\mathcal{L}^2}+(1+t)^2\|\widehat{u}\|^2_{L^2})+\lambda (1+t)^2\int_{\mathbb{R}^2\times B}\psi_\infty|\nabla_R(\frac{\psi-\psi_\infty}{\psi_\infty})|^2+(1+t)^2\int_{\mathbb{R}^2}|\xi|^2|\widehat{u}|^2d\xi\\
\leq 2(1+t)\lambda\|\psi-\psi_\infty\|^2_{\mathcal{L}^2}+2(1+t)\int_{S(t)}|\widehat{u}|^2d\xi.
\end{multline}
By virtue of (1.2), we get
\begin{align}\label{3.16}
\widehat{u}=e^{-t|\xi|^2}\widehat{u_0}+\int^t_0e^{-(t-s)|\xi|^2}i\xi \mathcal{F}(\mathbb{P}(u\otimes u)+\mathbb{P}\tau) ds,
\end{align}
where $\mathbb{P}$ stands for Leray's project operator. Using the fact that $|\widehat{f}|\leq \|f\|_{L^1}$, we have
\begin{align}\label{3.17}
|\widehat{u}|&\leq e^{-t|\xi|^2}|\widehat{u_0}|+|\xi|\int^t_0\|u\|^2_{L^2}ds+|\xi|t^{\frac{1}{2}}(\int^t_0|\widehat{\tau}|^2ds)^{\frac{1}{2}}\\
\nonumber&\leq C+|\xi|(\int^t_0\|u\|^3_{L^2}ds)^{\frac{1}{3}}t^{\frac{2}{3}}+|\xi|t^{\frac{1}{2}}(\int^t_0|\widehat{\tau}|^2ds)^{\frac{1}{2}}.
\end{align}

Using the system (1.2), we have
\begin{align}\label{3.18}
u=e^{-t\Delta}u_0+\int^t_0e^{-(t-s)\Delta}(\mathbb{P}(u\nabla u)+\mathbb{P}div \tau) ds.
\end{align}
Taking advantage of Lemma \ref{Lemma4}, we obtain
\begin{multline}
\|u\|_{L^2}\leq t^{-\frac{1}{4}}\|u_0\|_{L^{\frac{4}{3}}}+C\int^t_0(t-s)^{-\frac{1}{2}}\|u\nabla u\|_{L^1}+(t-s)^{-\frac{1}{2}}\|\tau\|_{L^2}ds\\
\leq \leq t^{-\frac{1}{4}}\|u_0\|_{L^{\frac{4}{3}}}+C\int^t_0(t-s)^{-\frac{1}{2}}\|u\|_{L^2}\|\nabla u\|_{L^2}+(t-s)^{-\frac{1}{2}}\|\tau\|_{L^2}ds.
\end{multline}
Note that $1+\frac{1}{3}=\frac{5}{6}+\frac{1}{2}$. Using the generalized Young inequality, we deduce that
\begin{align}\label{3.20}
(\int^t_0\|u(s)\|^3_{L^2}ds)^{\frac{1}{3}}&\leq t^{\frac{1}{12}} \|u_0\|_{L^{\frac{4}{3}}}+C\|s^{-\frac{1}{2}}1_{[0,t]}\|_{L^2_w}(\int^t_0(\|u(s)\|_{L^2}\|\nabla u(s)\|_{L^2}+\|\tau(s)\|_{L^2})^{\frac{6}{5}}ds)^{\frac{5}{6}}\\
\nonumber &\leq t^{\frac{1}{12}} \|u_0\|_{L^{\frac{4}{3}}}+C(\int^t_0\|u(s)\|^3_{L^2}ds)^{\frac{1}{3}}(\int^t_0\|\nabla u(s)\|^2_{L^2}ds)^{\frac{1}{2}}+(\int^t_0\|\tau(s)\|^{\frac{6}{5}}_{L^2}ds)^{\frac{5}{6}}\\
\nonumber &\leq  t^{\frac{1}{12}} \|u_0\|_{L^{\frac{4}{3}}}+aC(\int^t_0\|u(s)\|^3_{L^2}ds)^{\frac{1}{3}}+(\int^t_0\|\tau(s)\|^{\frac{6}{5}}_{L^2}ds)^{\frac{5}{6}},
\end{align}
where $a=(\|u_0\|^2_{L^2}+2K\|\psi_0-\psi_\infty\|^2_{\mathcal{L}^2})^{\frac{1}{2}}$. Applying Lemma \ref{Lemma3}, we verify that
\begin{align}\label{3.21}
(\int^t_0\|\tau(s)\|^{\frac{6}{5}}_{L^2}ds)^{\frac{5}{6}}&\leq \bigg\{\int^t_0\|\psi-\psi_\infty\|^{\frac{3}{5}}_{\mathcal{L}^2}[\int_{\mathbb{R}^2\times B}|\nabla_R(\frac{\psi-\psi_\infty}{\psi_\infty}|^2\psi_\infty dxdR]^{\frac{3}{5}}ds\bigg\}^{\frac{5}{6}}\\
\nonumber&\leq (\int^t_0\|\psi-\psi_\infty\|^{\frac{6}{7}}_{\mathcal{L}^2}ds)^{\frac{7}{12}}   \{\int^t_0[\int_{\mathbb{R}^2\times B}|\nabla_R(\frac{\psi-\psi_\infty}{\psi_\infty}|^2\psi_\infty dxdR]^2ds\}^{\frac{1}{4}}\\
\nonumber&\leq \|\psi_0-\psi_\infty\|_{\mathcal{L}^2}(\int^t_0\exp(-Ct))^{\frac{7}{12}}\leq C\|\psi_0-\psi_\infty\|_{\mathcal{L}^2}.
\end{align}
Plugging (\ref{3.21}) into (\ref{3.20}) yields that
\begin{align}
(\int^t_0\|u(s)\|^3_{L^2}ds)^{\frac{1}{3}}&\leq C(1+t)^{\frac{1}{12}} +aC(\int^t_0\|u(s)\|^3_{L^2}ds)^{\frac{1}{3}}.
\end{align}
 If $aC\leq \frac{1}{2}$, we then have
\begin{align}\label{3.23}
 (\int^t_0\|u(s)\|^3_{L^2}ds)^{\frac{1}{3}} &\leq C(1+t)^{\frac{1}{12}}.
\end{align}
Plugging (\ref{3.23}) into (\ref{3.17}) yields that
\begin{align}\label{3.24}
|\widehat{u}|\leq C+|\xi|(1+t)^{\frac{3}{4}}+|\xi|t^{\frac{1}{2}}(\int^t_0|\widehat{\tau}|^2ds)^{\frac{1}{2}},
\end{align}
which leads to
\begin{align}\label{3.25}
\int_{S(t)}|\widehat{u}|^2d\xi&\lm \int_{S(t)}d\xi+(1+t)^{\frac{3}{2}}\int_{S(t)}|\xi|^2d\xi+t\int_{S(t)}|\xi|^2(\int^t_0|\widehat{\tau}|^2ds)d\xi\\
\nonumber&\lm \int^{\sqrt{\frac{2}{1+t}}}_0rdr+(1+t)^{\frac{3}{2}}\int^{\sqrt{\frac{2}{1+t}}}_0r^{3}dr+\frac{2t}{1+t}\int^t_0\|\tau\|^2_{L^2}ds\\
\nonumber&\lm (1+t)^{-\frac{1}{2}}+\int^t_0\int_{\mathbb{R}^d\times B}\psi_\infty|\nabla_R(\frac{\psi-\psi_\infty}{\psi_\infty})|^2ds.
\end{align}
Plugging (\ref{3.25}) into (\ref{3.15}) and using the fact that $\|\psi-\psi_\infty\|_{\mathcal{L}^2}\lm \exp{(-Ct)}$ yield that
\begin{multline}
\frac{d}{dt}((1+t)^2\lambda\|\psi-\psi_\infty\|^2_{\mathcal{L}^2}+(1+t)^2\|\widehat{u}\|^2_{L^2})+\lambda (1+t)^2\int_{\mathbb{R}^2\times B}\psi_\infty|\nabla_R(\frac{\psi-\psi_\infty}{\psi_\infty})|^2\\
\leq C(1+t)^{\frac{1}{2}}+C(1+t)\int^t_0\int_{\mathbb{R}^2\times B}\psi_\infty|\nabla_R(\frac{\psi-\psi_\infty}{\psi_\infty})|^2ds.
\end{multline}
By taking $\lambda$ sufficiently large, we deduce that
\begin{align}
(1+t)^2\lambda\|\psi-\psi_\infty\|^2_{\mathcal{L}^2}+(1+t)^2\|u\|^2_{L^2}
\lm 1+ \int^t_0(1+t')^{\frac{1}{2}}dt'\lm (1+t)^{\frac{3}{2}},
\end{align}
which implies that
\begin{align}
\|u\|^2_{L^2}\lm (1+t)^{-\frac{1}{2}}.
\end{align}
From (\ref{3.16}) we have
\begin{align}
|\widehat{u}|&\leq e^{-t|\xi|^2}|\widehat{u_0}|+|\xi|\int^t_0\|u\|^2_{L^2}ds+|\xi|t^{\frac{1}{2}}(\int^t_0|\widehat{\tau}|^2ds)^{\frac{1}{2}}\\
\nonumber&\leq
\|u_0\|_{L^1}+C|\xi|\int^t_0(1+s)^{-\frac{1}{2}}ds+|\xi|t^{\frac{1}{2}}(\int^t_0|\widehat{\tau}|^2 ds)^{\frac{1}{2}}\\
\nonumber&\leq \|u_0\|_{L^1}+C|\xi|\int^t_0(1+s)^{-\frac{1}{2}}ds+|\xi|t^{\frac{1}{2}}(\int^t_0|\widehat{\tau}|^2 ds)^{\frac{1}{2}}\\
\nonumber&=\|u_0\|_{L^1}+C|\xi|\sqrt{1+t}+|\xi|t^{\frac{1}{2}}(\int^t_0|\widehat{\tau}|^2 ds)^{\frac{1}{2}},
\end{align}
which leads to
\begin{align}\label{3.30}
\int_{S(t)}|\widehat{u}|^2d\xi&\lm \int_{S(t)}d\xi+(1+t)\int_{S(t)}|\xi|^2d\xi+t\int_{S(t)}|\xi|^2(\int^t_0|\widehat{\tau}|^2ds)d\xi\\
\nonumber&\lm \int^{\sqrt{\frac{2}{1+t}}}_0rdr+(1+t)\int^{\sqrt{\frac{2}{1+t}}}_0r^{3}dr+\frac{2t}{1+t}\int^t_0\|\tau\|^2_{L^2}ds\\
\nonumber&\lm (1+t)^{-1}+\int^t_0\int_{\mathbb{R}^2\times B}\psi_\infty|\nabla_R(\frac{\psi-\psi_\infty}{\psi_\infty})|^2ds.
\end{align}
Plugging (\ref{3.30}) into (\ref{3.15}) yields that
\begin{multline}
\frac{d}{dt}((1+t)^2\lambda\|\psi-\psi_\infty\|^2_{\mathcal{L}^2}+(1+t)^2\|\widehat{u}\|^2_{L^2})+\lambda (1+t)^2\int_{\mathbb{R}^2\times B}\psi_\infty|\nabla_R(\frac{\psi-\psi_\infty}{\psi_\infty})|^2\\
\leq C+C(1+t)\int^t_0\int_{\mathbb{R}^2\times B}\psi_\infty|\nabla_R(\frac{\psi-\psi_\infty}{\psi_\infty})|^2ds.
\end{multline}
By taking $\lambda$ sufficiently large, we get
\begin{align}
(1+t)^2\lambda\|\psi-\psi_\infty\|^2_{\mathcal{L}^2}+(1+t)^2\|u\|^2_{L^2}
\lm 1+ \int^t_0dt'\lm (1+t),
\end{align}
which implies that
\begin{align}
\|u\|^2_{L^2}\lm (1+t)^{-1}.
\end{align}
\end{proof}
\end{theo}

\subsection{The $L^2$ decay with large data}
\begin{theo}\label{th3}
 Suppose that $p\in[1,\infty]$ and $pk>1$. Let $(u,\psi)$ be a weak solution of (1.2) with the initial data $u_0\in L^2\cap \dot{B}^0_{1,1}$ and $\psi_0$ satisfies $\psi_0-\psi_\infty\in L^2_x(\mathcal{L}^2)\cap L^1_x(\mathcal{L}^p)$ and $\int_B\psi_0=1$ $a.e.$ in $x$.
Then there exists a constants such that
\begin{align}
\int_{\mathbb{R}^2\times B}\frac{|\psi-\psi_\infty|^2}{\psi_\infty}dxdR\leq C\exp{(-Ct)},
\end{align}
\begin{align}
\|u\|_{L^2(\mathbb{R}^2)}\leq C(1+t)^{-\frac{1}{2}}.
\end{align}
\end{theo}

In order to prove the above theorem, we need to use the Littlewood-Paley decomposition and some basic lemma for the homogeneous Besov space. (see \cite{Bahouri2011} for more details)

Let $\mathcal{C}$ be the annulus $\{\xi\in\mathbb{R}^{d}\big|\frac{3}{4}\leq|\xi|\leq\frac{8}{3}\}.$ There exists radial function $\varphi$, valued in the interval $[0,1]$, such that
\begin{align}
\forall\xi\in\mathbb{R}^{d}\backslash\{0\},~\sum_{j\in\mathbb{Z}}\varphi(2^{-j}\xi)=1,
\end{align}
\begin{align}
|j-j'|\geq2\Rightarrow Supp ~\varphi(2^{-j}\xi)\cap Supp ~\varphi(2^{-j'}\xi)=\emptyset.
\end{align}
 The homogeneous dyadic blocks $\dot{\Delta}_{j}$ are defined by
\begin{align} \dot{\Delta}_{j}u=\varphi(2^{-j}D)u=2^{jd}\int_{\mathbb{R}^{d}}h(2^{j}y)u(x-y)dy,
\end{align}
\begin{align}
\dot{S}_{j}u=\chi(2^{-j}D)u=\int_{\mathbb{R}^{d}}\widetilde{h}(2^{j}y)u(x-y)dy.
\end{align}

   The homogeneous Besov space is denoted by $\dot{B}^{s}_{p,r}$, that is
$$\dot{B}^{s}_{p,r}=\big\{u\in S'_{h}\big{|}\|u\|_{\dot{B}^{s}_{p,r}}=\|2^{js}\|\dot{\Delta}_{j}u\|_{L^{p}_x}\|_{l^r}<\infty\big\},$$

The following lemmas will be useful to obtain the estimates for the solutions of the Navier-Stokes equations.
\begin{lemm}\label{L1}
\cite{Bahouri2011} Let $\mathcal{C}$ be an annulus and $B$ a ball. A constant $C$ exists such that
for any nonnegative integer $k$, any couple $(p, q)$ in $[1,\infty]^{2}$ with $q \geq p\geq 1$, and
any function $u$ of $L^{p}$, we have
$$Supp ~\widehat{u} \subseteq \lambda B \Rightarrow \|D^{k}u\|_{L^{q}}\triangleq \sup_{|\alpha|\leq k} \|\partial ^{\alpha}u\|_{L^{q}}\leq C^{k+1}\lambda^{k+d(\frac{1}{p}-\frac{1}{q})}\|u\|_{L^{p}},$$
$$Supp ~\widehat{u} \subseteq \lambda\mathcal{C} \Rightarrow C^{-k-1}\lambda^{k}\|u\|_{L^{p}} \leq\|D^{k}u\|_{L^{p}}\leq C^{k+1}\lambda^{k}\|u\|_{L^{p}}.$$
$$Supp ~\widehat{u} \subseteq \lambda\mathcal{C} \Rightarrow  \|e^{t\triangle}u\|_{L^{p}}\leq Ce^{-ct\lambda^{2}}\|u\|_{L^{p}}.$$
\end{lemm}

\begin{lemm}\label{L2}
\cite{Bahouri2011} For any positive $s$, we have
$$\sup_{t>0}\sum_{j}t^s2^{2js}e^{-ct2^{2j}}<\infty.$$
\end{lemm}

\begin{lemm}\label{Lemma5}
Suppose that $u_0\in L^2\cap \dot{B}^0_{1,1}$ and $\tau\in L^2_T(L^2)\cap L^\infty_T(L^1)$. If $u\in L^\infty_T(L^2)\cap L^2_T(\dot{H}^1)$ is the solution of
   \begin{align}\label{3.40}
\left\{
\begin{array}{ll}
u_t+(u\cdot\nabla)u-\Delta{u}+\nabla{P}=div~\tau, ~~~~~~~div~u=0,\\[1ex]
u|_{t=0}=u_0,
\end{array}
\right.
\end{align}
then we have
$$\sup_{t\in[0,T]}\|u\|_{\dot{B}^0_{1,1}}\leq \|u_0\|_{\dot{B}^0_{1,1}}+C\sqrt{T}(\|u_0\|^2_{L^2}+\int^T_0\|\tau\|^2_{L^2}ds+\|\tau\|_{L^\infty_T(L^1)}).$$
\begin{proof}
By the standard energy estimate, we have
\begin{align}
\|u\|^2_{L^2}+\frac{1}{2}\int^T_0\|\nabla u\|^2_{L^2}ds\leq \|u_0\|^2_{L^2}+\int^T_0\|\tau\|^2_{L^2}ds.
\end{align}
By virtue of (\ref{3.40}), we can write that
\begin{align}
u=e^{-t\Delta}u_0+\int^t_0e^{-(t-s)\Delta}(\mathbb{P}div (u\otimes u+ \tau) ds.
\end{align}
Applying $\dot{\Delta}_j$ to the above equation and taking the $L^1-$norm, we deduce from Lemma \ref{L1} that
\begin{multline}
\|\dot{\Delta}_ju\|_{L^1}\leq \|\dot{\Delta}_ju_0\|_{L^1}+\int^t_0e^{-c(t-s)2^{2j}}2^j(\|\dot{\Delta}_j(u\otimes u)\|_{L^1}+\|\dot{\Delta}_j\tau\|_{L^1})ds\\
\leq \|\dot{\Delta}_ju_0\|_{L^1}+\int^t_0e^{-c(t-s)2^{2j}}2^j(\|(u\otimes u)\|_{L^1}+\|\tau\|_{L^1})ds,
\end{multline}
which leads to
\begin{multline}
\sum_{j}\|\dot{\Delta}_ju\|_{L^1}\leq \|u_0\|_{\dot{B}^0_{1,1}}+\int^t_0(t-s)^{-\frac{1}{2}}\sum_j(t-s)^{\frac{1}{2}}e^{-c(t-s)2^{2j}}2^j(\|(u\otimes u)\|_{L^1}+\|\tau\|_{L^1})ds\\
\leq \|u_0\|_{\dot{B}^0_{1,1}}+\int^t_0(t-s)^{-\frac{1}{2}}ds\sup_{t-s>0}\sum_j(t-s)^{\frac{1}{2}}e^{-c(t-s)2^{2j}}2^j(\|u\|^2_{L^\infty_T(L^2)}+\|\tau\|_{L^\infty_T(L^1)})\\
\leq \|u_0\|_{\dot{B}^0_{1,1}}+C\sqrt{T}(\|u_0\|^2_{L^2}+\int^T_0\|\tau\|^2_{L^2}ds+\|\tau\|_{L^\infty_T(L^1)}).
\end{multline}
\end{proof}
\end{lemm}

\textit{Proof of Theorem 3.2:} Multiplying $p|\frac{\psi-\psi_\infty}{\psi_\infty}|^{p-2}\frac{\psi-\psi_\infty}{\psi_\infty}$ by both sides of (\ref{3.4}) and integrating over $B$ with $R$, we obtain
\begin{multline}\label{3.45}
\frac{d}{dt}\int_B|\frac{\psi-\psi_\infty}{\psi_\infty}|^p\psi_\infty+u\cdot\nabla_x\int_B|\frac{\psi-\psi_\infty}{\psi_\infty}|^p\psi_\infty+\frac{4(p-1)}{p}\int_B\psi_\infty|\nabla_R(\frac{\psi-\psi_\infty}{\psi_\infty})^{\frac{p}{2}}|^2
\\
=\int_B\sigma(u)R(\psi-\psi_\infty)\nabla_R(\frac{\psi-\psi_\infty}{\psi_\infty})^p.
\end{multline}
Using integration by parts and (\ref{3.3}), we see that
\begin{align}\label{3.5}
\frac{d}{dt}\int_B|\frac{\psi-\psi_\infty}{\psi_\infty}|^p\psi_\infty+u\cdot\nabla_x\int_B|\frac{\psi-\psi_\infty}{\psi_\infty}|^p\psi_\infty+\frac{4(p-1)}{p}\int_B\psi_\infty|\nabla_R(\frac{\psi-\psi_\infty}{\psi_\infty})^{\frac{p}{2}}|^2=0,
\end{align}
which leads to
\begin{align*}
\frac{d}{dt}\int_B|\frac{\psi-\psi_\infty}{\psi_\infty}|^p\psi_\infty+u\cdot\nabla_x\int_B|\frac{\psi-\psi_\infty}{\psi_\infty}|^p\psi_\infty\leq 0.
\end{align*}
Since $div~u=0$, it follows that
$$\|\psi-\psi_\infty\|_{L^1_x\mathcal{L}^p}\leq \|\psi_0-\psi_\infty\|_{L^1_x\mathcal{L}^p}.$$
Taking advantage of H\"{o}lder's inequality and using the fact that $pk>1$, we have
\begin{align}
|\tau|\leq \int_{B}\frac{|\psi-\psi_\infty|}{1-|R|}dR\leq \int_{B}\frac{(\psi_\infty)^{\frac{p}{p+1}}}{(1-|R|)}\frac{|\psi-\psi_\infty|}{(\psi_\infty)^\frac{p}{p+1}}dR\leq C (\int_B|\frac{\psi-\psi_\infty}{\psi_\infty}|^p\psi_\infty dR)^{\frac{1}{p}},
\end{align}
which leads to
\begin{align}
\|\tau\|_{L^1}\leq \|\psi-\psi_\infty\|_{L^1_x\mathcal{L}^p}\leq \|\psi_0-\psi_\infty\|_{L^1_x\mathcal{L}^p}.
\end{align}
By virtue of Lemma \ref{Lemma5}, we deduce that
\[\sup_{t\in[0,T]}\|u\|_{\dot{B}^0_{1,1}}\leq C_T,\]
for any $T<\infty$. By virtue of Theorem \ref{th1}, we see that
\[\|u\|_{L^2}\leq C\ln^{-1}(e+t), quad \|\psi-\psi_\infty\|_{\mathcal{L}^2}\leq C\exp{-Ct},\]
which implies that for any $\ep>0$ there exists $T_0$ such that
\[\|u(T_0)\|_{L^2}+\|\psi(T_0)-\psi_\infty\|_{\mathcal{L}^2}<\ep.\]
Since $\dot{B}^0_{1,1}\hookrightarrow L^1$. it follows that $\|u(T_0)\|_{L^1}\leq C_{T_0}$.
Let $(u(T_0),\psi(T_0))$ be the initial data. Applying Theorem \ref{th2}, we complete the proof.

\begin{rema}
By Theorem \ref{th3} together with Theorem \ref{th1}, we see that the conjecture proposed by M. Schonbek in \cite{Schonbek} holds true for all $d\geq 2$.
In \cite{Schonbek1991}, M. Schonbek showed that $(1+t)^{-\frac{d}{4}},\, d\geq 2$ is the optimal $L^2$ decay rate for the Navier-Stokes equations with $u_0\in L^1$. Note that
if $\psi$ is independent on $x$, then $div~\tau=0$. Then, the co-rotation FENE model is reduced to the Navier-Stokes equations. Thus, the $L^2$ decay rate for the co-rotation FENE model which we obtained in Theorem \ref{th1} and Theorem \ref{th3} is sharp for all $d\geq 2$.
\end{rema}
\smallskip
\noindent\textbf{Acknowledgments} This work was
partially supported by the National Natural Science Foundation of China (No.11671407 and No.11701586), the Macao Science and Technology Development Fund (No. 098/2013/A3), and Guangdong Province of China Special Support Program (No. 8-2015), 
and the key project of the Natural Science Foundation of Guangdong province (No. 2016A030311004).


\phantomsection
\addcontentsline{toc}{section}{\refname}
\bibliographystyle{abbrv} 
\bibliography{Feneref}

\end{document}